\documentclass{article}

\usepackage{amsmath, amsthm, amssymb, verbatim, enumerate, tikz-cd, graphicx, float}

\title{The Spectrum of the Quaquaversal Operator is Real}
\author{Josiah Sugarman}
\date{}

\title{The Spectrum of the Quaquaversal Operator is Real}
\author{Josiah Sugarman}
\date{}

\DeclareMathOperator{\SO3}{SO(3)}
\DeclareMathOperator{\SOn}{SO(n)}

\DeclareMathOperator{\tr}{tr}

\newtheorem{theorem}{Theorem}

\newtheorem*{corollary}{Corollary}
\newtheorem{remark}{Remark}
\newtheorem{lemma}{Lemma}
\begin{document}
\maketitle{}
\begin{abstract}
The Hecke operator associated the Quaquaversal tiling introduced by John Conway and Charles Radin is shown to have a real spectrum--- resolving a conjecture of Draco, Sadun, and Van Wieren. 
\end{abstract}
\section{Background and Introduction}\label{sec:Intro}
Suppose a {\em mother tile} is tiled by smaller copies of itself--{\em daughter tiles}--as in the following example:
\begin{figure}[H]
\centering
\includegraphics[height = 2cm]{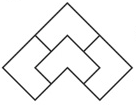}
\caption{L-shaped mother tile (Chair tile) tiled by its four  daughters}
\end{figure}
A tiling of the mother tile by daughter tiles like this gives rise to an infinite tiling as follows:

Choose a daughter tile. The mother tile will play the role of this daughter while tiling the {\em grandmother tile} the grandmother will play the same role when tiling the {\em great grandmother tile} and the process will continue indefinitely, tiling all the ancestors and ancient ancestors as in the following figure:
\begin{figure}[H]
\centering
\includegraphics[]{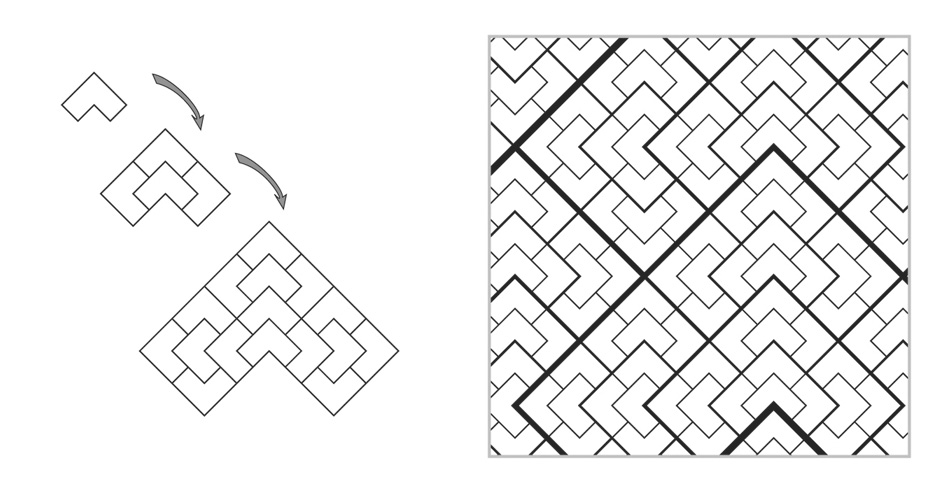}
\caption{The Chair tiling of the plane}
\end{figure}
We call such a tiling {\em hierarchical}. More precisely, a hierarchical tiling is a tile, $T\subseteq \mathbb{R}^N$, together with $k$ affine maps, $\{g_0,\ldots, g_{k-1}\}$, of the form $c_i + \lambda_i O_i$ with $0 < \lambda_i < 1$, $c_i\in\mathbb{R}^n$, and $O_i\in \SOn$ such that $T=\cup_{i=0}^{k-1} g_i T$ with $g_iT$ disjoint except perhaps on their boundaries. We call $T$ a mother tile and the $g_iT$ are $T$'s daughter tiles.

Without loss of generality, let $T$ play the role of $g_0T$ when tiling it's mother, $g_0^{-1}T$. We may construct an infinite tiling as follows:

Since
\[T = \cup_{i=0}^{k-1} g_i T = \cup_{i,j=0}^{k-1} g_ig_j T =\cdots =\cup_{w\in G_N} w T \]

where $G_N$ is the set of words in $g_0, g_1,\ldots, g_{k-1}$ of length $N$ we have:

\[g_0^{-N}T = \cup_{w\in G_N} g_0^{-N}wT. \]

Since $T \subset g_0^{-1}T\subset \cdots \subset g_0^{-N}T\subset\cdots$ we have a tiling of $\cup_{i=0}^{\infty}g_0^{-i}T$.
As $T$ is the mother of the $g_iT$, it makes sense to call tiles of the form $g_0^{-N}wT$ with $w\in G_N$ the (N-1)\textsuperscript{st} cousins of $T$. Naturally the `sibling tiles' of $T$, $g_0^{-1}g_iT$, are $T$'s `0\textsuperscript{th} cousins' and $T$ is $T$'s -1\textsuperscript{st} cousin. We can summarize the above discussion by saying the cousins of $T$ tile $\cup_{i=0}^{\infty}g_0^{-i}T$.

Notice that if $w = g_{i_1}^{a_1}g_{i_2}^{a_2}\cdots g_{i_N}^{a_N}$ with $g_{i_k} = \lambda O_{i_k} + c_{i_k}$ then the `rotational part' of $w$ is given by $O_{i_1}^{a_1}O_{i_2}^{a_2}\cdots O_{i_N}^{a_N}$. By identifying the orientation of $T$ with the identity element of $SO(n)$ we now see that the orientations of the cousin tiles are given by $O_0^{-(a_1+\cdots+a_N)} O_{i_1}^{a_1}O_{i_2}^{a_2}\cdots O_{i_N}^{a_N}$. When $O_0$ is the identity (when the daughter tile has the same orientation as the mother tile) this simplifies to $O_{i_1}^{a_1}O_{i_2}^{a_2}\cdots O_{i_N}^{a_N}$.


In 1998, Charles Radin and John Conway introduced a three dimensional hierarchical tiling, {\em the Quaquaversal tiling} \cite{conway1998quaquaversal}. This tiling exhibits ``statistical rotational symmetry in the infinite volume limit'' meaning the distribution of orientations of the tiles chosen uniformly at random from a sphere approach uniformity as the radius of the sphere approaches infinity. 

The Quaquaversal tiling is particularly interesting because this distribution approaches uniformity at a much faster rate than what is possible for two dimensional hierarchical tilings. Indeed, $SO(2)$ is abelian so the number of distinct words in any generating set grows as a polynomial with respect to the word-length. But, $\SO3$ is highly nonabelian and contains many generating sets of exponential growth. This gives three dimensional tilings the opportunity for tiles to rapidly equidistribution.


The Quaquaversal tiling is a hierarchical tiling whose mother tile is a right angled triangular prism with depth 1, height 1 and length $\sqrt{3}$ and partitioned as indicated:

\begin{figure}[H]
\includegraphics[height = 4cm]{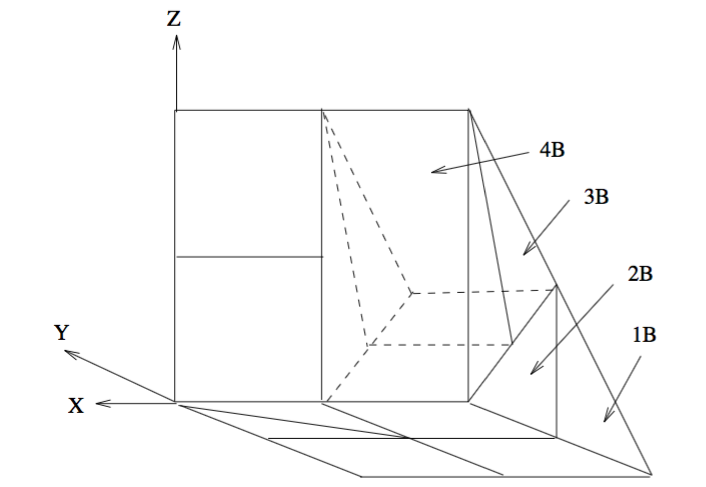}
\includegraphics[height = 4cm]{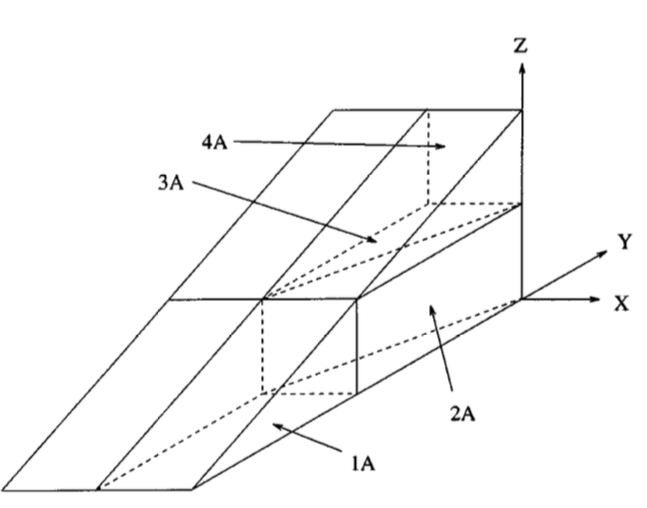}
\caption{The Quaquaversal Tiling \protect\cite{conway1998quaquaversal} \protect\cite{draco2000growth}}
\end{figure}

After identifying the orientation of the mother tile with the identity rotation, the orientations of the eight daughter tiles are given by three identities, $S^2T^3$, $T^4$, $T^4S^2$, $S$, and $ST^3$ where $S$ and $T$ are rotations about orthogonal axes by $\pi/2$ and $\pi/3$ respectively. Taking a formal sum of these rotations and dividing by $8$ gives us an element of the group ring $\mathbb{C}[\SO3]$ that acts as a sort of `generating function' for the distribution of the cousins' orientations, in the sense that if you raise this formal sum to the $N$ then you get the distribution of the orientations of the $(N-1)^{\text{st}}$ cousins. Such a formal sum gives rise to an operator on $L^2(\SO3)$.

By showing the eigenvalue $1$ occurs with multiplicity equal to $1$ for this operator, Conway and Radin showed that this tiling exhibits asymptotic rotational symmetry in the above sense. By finding eigenvalues larger than $0.9938$ Draco, Sadun, and Van Wieren\cite{draco2000growth} showed that the rate with which this tiling approaches uniformity is quite slow---mentioning that a galactic scale compound of a chemical with this structure would still exhibit noticeable anisotropy; By showing there is a positive gap separating $1$, the largest eigenvalue, from all other eigenvalues, Bourgain and Gamburd\cite{bourgain2008spectral} showed that this tiling approaches uniformity at an exponentially fast rate.

In the course of their numerical experiments, Draco, Sadun, and Van Wieren\cite{draco2000growth} observed and conjectured that this operator has a real spectrum. In section 3, we state their conjecture filling in the gaps necessary to define the operator. In section 4 we prove this conjecture. The main ingredient of the proof is a well chosen partition of the representation. With respect to this partition, the operator is lower block triangular with Hermitian blocks along the main diagonal. In the fifth section the blocks are described in more detail and in the last section the spectrum of three of the four blocks along the main diagonal are analyzed.



\subsection{Draco, Sadun, and Van Wieren Conjecture}\label{sec:SDWConj}
Combining the Peter-Weyl theorem with the classification of irreducible representations of $\SO3$ \cite{bump2004lie}\cite{knapp2016representation}\cite{weyl2016classical} we get:
\[L^2(\SO3) = \oplus_{2k + 1} H_{2k+1}^{2k+1}\]
where $k$ ranges over the positive integers and the $H_{2k+1}$ are the unique $2k+1$ dimensional representations of dimension $2k+1$. This direct sum is orthogonal with respect to an $\SO3$ invariant inner product so that the linear operator $v\mapsto gv$ is unitary. Furthermore, the induced operators on $L^2(\SO3)$ coming from elements of the group ring, $\mathbb{C}[\SO3]$, leave the $H_{2k+1}$ invariant. So an element $z\in\mathbb{C}[\SO3]$ gives us an infinite family of linear operators, \[\pi_{2k+1}(z):H_{2k+1}\rightarrow H_{2k+1}\]
\paragraph{Draco, Sadun, Van Wieren Conjecture:}
The quaquaversal tiling is a hierarchical tiling with associated element of the group ring given by:
\[z = 1/8(1 + 1 + 1 + S^2T^3 + T^4 + T^4S^2 + S + ST^3)\]
where $S = R_a^{\pi/2}$ and $T=R_b^{\pi/3}$ are rotations about orthogonal vectors, $a$, and $b$, by $\pi/2$ and $\pi/3$, respectively.
Draco, Sadun, and Van Wieren conjectured that the associated operator
$\pi(z):L^2(\SO3)\rightarrow L^2(\SO3)$, or equivalently the infinite family of finite dimensional operators
\[\pi_{2k+1}(z):\mathcal{H}_{2k+1}\rightarrow\mathcal{H}_{2k+1}\]
have a real spectrum.

\subsection{Notation}
From now on we shall work in an arbitrary, yet fixed, finite dimensional, unitary, irreducible representation of $\SO3$ of dimension $2k+1$; we will let \mbox{$R_b^\theta$} denote the unitary representation of dimension $2k+1$ of the rotation about the $v$ axis by $\theta$; and we will write $\hat{z}$ for $\pi_{2k+1}(z)$.

We shall use the following nonstandard notation for block partitioned matrices. Let $M$ be a $(2k+1) \times (2k+1)$ matrix partitioned with respect to $\Pi = V_1,\ldots, V_d$. There are unique projections $\{P_i\}$ which are the identity on $V_i$ and vanish on the $V_j$ for $j\neq i$.
When we write
\[M=
\begin{pmatrix}
A_{1,1}& \cdots &A_{1,j} & \cdots & A_{1,d}\\
\vdots& \ddots &\vdots & \ddots &\vdots\\
A_{i,1}& \cdots &A_{i,j} & \cdots & A_{i,d}\\
\vdots& \ddots &\vdots & \ddots &\vdots\\
A_{d,1}& \cdots &A_{d,j} & \cdots & A_{d,d}\\
\end{pmatrix}
\]
With the $A_{i,j}$'s being $(2k+1)\times (2k+1)$ matrices we mean that $M|_{V_i}^{V_j} = P_j A P_i$. Additionally, we require that the $A_{i,j}$ map $V_i$ to $V_j$ so that $M|_{V_i}^{V_j} = P_j A P_i = A|_{V_i}$. That is, for $v\in V_i$, $P_jMv = Av$.
\subsection{Acknowledgements}
I'd like to thank Alexander Gamburd for suggesting this problem and many helpful conversation and Lorenzo Sadun for noticing a small error in an earlier version.
\section{Proof}\label{sec:main_thm}
\begin{theorem}\label{thm:main_thm}
Let $a$ and $b$ be orthogonal vectors. Then $R_a^\theta + R_a^\theta R_b^\pi$ is of the form:
\[\begin{bmatrix}
2R_a^{\theta+} & 0& 0& 0 \\
2R_a^{\theta-} & 0& 0& 0 \\
0& 0& 2R_a^{\theta+}& 0 \\
0& 0& 2R_a^{\theta-} & 0
\end{bmatrix}
\]
When partitioned with respect to 
\[\Pi_{a,b} = \Lambda_{1,1},\Lambda_{1, -1}, \Lambda_{-1, 1}, \Lambda_{-1,-1}\]
with, 
\[\Lambda_{\alpha, \beta}:=\{v| R_a^\pi v =\alpha v\text{ and } R_b^\pi =\beta v\}\]
where $R_a^{\theta+}$, the {\em Hermitian part} of $R_a^\theta$, is equal to $\frac{R_a^\theta+R_a^{-\theta}}{2}$ and $R_a^{\theta-}$, the {\em skew-Hermitian part} of $R_a^\theta$, is equal to $\frac{R_a^\theta-R_a^{-\theta}}{2}$.
\end{theorem}
\begin{corollary}
The spectrum of the Quaquaversal Operator, $z = 1/8(S + ST^3 + S^2T^3 + T^4 +T^4S^2 + 3) $, is real.
\end{corollary}
\begin{proof}
By theorem \ref{thm:main_thm}, $S+ST^3$ can be block triangulated as
\[\begin{bmatrix}
2R_a^{\pi/2+} & 0& 0& 0 \\
2R_a^{\pi/2-} & 0& 0& 0 \\
0& 0& 2R_a^{\pi/2+}& 0 \\
0& 0& 2R_a^{\pi/2-} & 0
\end{bmatrix}
\]
with respect to $\Pi_{a, b}$. Similarly, by theorem \ref{thm:main_thm}, $T^4+T^4S^2$ can be block triangulated as

\[\begin{bmatrix}
2R_b^{4\pi/3+} & 0& 0& 0 \\
2R_b^{4\pi/3-} & 0& 0& 0 \\
0& 0& 2R_b^{4\pi/3+} & 0 \\
0& 0& 2R_b^{4\pi/3-} & 0
\end{bmatrix}
\]
with respect to $\Pi_{b, a}$. And this is 
\[\begin{bmatrix}
2R_b^{4\pi/3+} & 0& 0& 0 \\
0 & 2R_b^{4\pi/3+}& 0& 0 \\
2R_b^{4\pi/3-}& 0& 0 & 0 \\
0& 2R_b^{4\pi/3-}& 0 & 0
\end{bmatrix}
\]
with respect to $\Pi_{a, b}$, since $\Pi_{a, b}$ is just $\Pi_{b, a}$ with $\Lambda_{-1, 1}$ and $\Lambda_{1, -1}$ switched (therefore, the blocks of the second and third `rows' and `columns' are swapped). Adding these we get that $S+ ST^3 + T^4 + T^4S^2$ is

\begin{equation}
\begin{bmatrix}
2R_a^{\pi/2+} + 2R_b^{4\pi/3+}& 0& 0& 0 \\
2R_a^{\pi/2-} & 2R_b^{4\pi/3+}& 0& 0 \\
2R_b^{4\pi/3-}& 0& 2R_a^{\pi/2+}& 0 \\
0& 2R_b^{4\pi/3-}& 2R_a^{\pi/2-} & 0
\end{bmatrix}
\end{equation}
with respect to $\Pi_{a,b}$. On $\Lambda_{\alpha, \beta}$, $S^2T^3 = \alpha\beta I$. So, $S + ST^3 + T^4 + T^4S^2 + S^2T^3$ is 
\[
\begin{bmatrix}
2R_a^{\pi/2+} + 2R_b^{4\pi/3+}+ I& 0& 0& 0 \\
2R_a^{\pi/2-} & 2R_b^{4\pi/3+}-I& 0& 0 \\
2R_b^{4\pi/3-}& 0& 2R_a^{\pi/2+}-I& 0 \\
0& 2R_b^{4\pi/3-}& 2R_a^{\pi/2-} & I
\end{bmatrix}
\]
As the blocks along the main diagonal are Hermitian and the matrix is lower block triangular, this means the spectrum of $S+ ST^3 + T^4 + T^4S^2 + S^2T^3$ is real. Since $1/8(S+ ST^3 + T^4 + T^4S^2 + S^2T^3 + 3I)$ only shifts and rescales the eigenvalues, it too must have a real spectrum, as needed.
\end{proof}
\begin{remark}\label{rem:BlockStruc}
Since we have block triangulated $S+ ST^3 + T^4 + T^4S^2 + S^2T^3$, we get a block triangulation of $z  = 1/8(S+ ST^3 + T^4 + T^4S^2 + S^2T^3 + 3I)$, the Quaquaversal Operator. It is block triangulated with respect to $\Pi_{a, b}$ as follows :
\[1/8\begin{bmatrix}
2R_a^{\pi/2+} + 2R_b^{4\pi/3+}+4I& 0& 0& 0 \\
2R_a^{\pi/2-}& 2R_b^{4\pi/3+}+2I& 0& 0 \\
2R_b^{4\pi/3-}& 0& 2R_a^{\pi/2+}+2I& 0 \\
0& 2R_b^{4\pi/3-}& 2R_a^{\pi/2-} & 4I
\end{bmatrix}
\]
When $(R_a^{\pi/2})^2 = -I$,  $R_a^{\pi/2} = -R_a^{-\pi/2}$. So $R_a^{\pi/2+} = \frac{R_a^{\pi/2}+R_a^{-\pi/2}}{2}= \frac{R_a^{\pi/2}-R_a^{\pi/2}}{2} = 0$ when the domain is restricted to $\Lambda_{-1, 1}$. When $(R_a^{\pi/2})^2 = I$,  $R_a^{\pi/2} = R_a^{-\pi/2}$. So $R_a^{\pi/2-} = \frac{R_a^{\pi/2}-R_a^{-\pi/2}}{2}= \frac{R_a^{\pi/2}-R_a^{\pi/2}}{2} = 0$ when the domain is restricted to $\Lambda_{1, 1}$. Furthermore, $R_b^{\pi} = I$ on $\lambda_{\alpha, 1}$ so $R_b^{4\pi} = R_b^{\pi}$ on $\lambda_{\alpha, 1}$ and $R_b^{\pi} = -I$ on $\lambda_{\alpha, -1}$ so $R_b^{4\pi} = -R_b^{\pi}$ on $\lambda_{\alpha, -1}$. Putting this together, the block triangularization simplifies to

\[1/8\begin{bmatrix}
2R_a^{\pi/2+} + 2R_b^{\pi/3+}+4I& 0& 0& 0 \\
0& -2R_b^{\pi/3+}+2I& 0& 0 \\
2R_b^{\pi/3-}& 0& 2I& 0 \\
0& -2R_b^{\pi/3-}& 2R_a^{\pi/2-} & 4I
\end{bmatrix}
\]
\end{remark}
\begin{proof}[Proof of Theorem \ref{thm:main_thm}]
Consider the subspaces $\Lambda_{\alpha, \beta}$ and the partition, $\Pi_{a, b}$, as in the statement of this theorem.
\begin{lemma}\label{lem:eigen}
If $A$ commutes with $B$ then $A$ preserves $B$'s eigenspace. That is, $AB = BA$ implies $Bv = \lambda v \Rightarrow BAv = \lambda Av$.

If $A$ anticommutes with $B$ then $A$ maps $B$'s $\lambda$ eigenspace to $B$'s $-\lambda$ eigenspace. That is, $AB = -BA$ implies $Bv =\lambda v\Rightarrow BAv = -\lambda Av $
\end{lemma}
\begin{proof}[proof of Lemma]
Suppose $AB = BA$ and $Bv = \lambda v$ then $BAv = ABv = A\lambda v = \lambda A v$ as needed.

Suppose $AB = -BA$ and $Bv = \lambda v$ then $BAv = -ABv = -A\lambda v = -\lambda A v$ as needed.
\end{proof}
Decompose $R_a^\theta$ into it's Hermitian and skew-Hermitian parts:
\[R_a^\theta = \underbrace{\frac{R_a^\theta+R_a^{-\theta}}{2}}_{R_a^{\theta+}}+\underbrace{\frac{R_a^\theta-R_a^{-\theta}}{2}}_{R_a^{\theta-}}\]

Note that $R_a^\pi R_a^{\theta} = R_a^{\theta}R_a^\pi$ and $R_b^\pi R_a^{\theta} = R_a^{-\theta}R_b^\pi$. Therefore $R_a^{\theta+}$ commutes with both $R_a^\pi$ and $R_b^\pi$. Therefore, by lemma {\ref{lem:eigen}}, $R_a^{\theta+}$ preserves both of their eigenspaces and therefore:
\[R_a^{\theta+}:\Lambda_{\alpha,\beta}\rightarrow\Lambda_{\alpha, \beta}\]

Similarly, since $R_a^{\theta-}$ commutes with $R_a^\pi$ and anticommutes with $R_b^\pi$, by lemma {\ref{lem:eigen}} we have:
\[R_a^{\theta-}:\Lambda_{\alpha,\beta}\rightarrow\Lambda_{\alpha,-\beta}\]

 We can now block partition $R_a^\theta = R_a^{\theta+}+R_a^{\theta-}$ with respect to $\Pi_{a, b}$
\[R_a^\theta = \left[\begin{matrix}
R_a^{\theta+}& R_a^{\theta-}& 0& 0 \\
R_a^{\theta-}& R_a^{\theta+}& 0& 0\\
0& 0& R_a^{\theta+}&  R_a^{\theta-}\\
0& 0& R_a^{\theta-} & R_a^{\theta+}
\end{matrix}\right]
\]

Now consider $R_a^\theta + R_a^\theta R_b^\pi = R_a^\theta(I + R_b^\pi)$. When we restrict the domain to either $\Lambda_{\alpha, -1}$, $R_b^\pi = -I$, therefore $R_a^\theta(I + R_b^\pi) = R_a^\theta(I - I) = 0$. When we restrict the domain to either $\Lambda_{\alpha, 1}$, $R_b^\pi = I$, therefore $R_a^\theta(I + R_b^\pi) = R_a^\theta(I + I) = 2R_a^\theta$. So,

\[R_a^\theta + R_a^\theta R_b^\pi = \left[\begin{matrix}
2R_a^{\theta+}& 0& 0& 0 \\
2R_a^{\theta-}& 0& 0& 0\\
0& 0& 2R_a^{\theta+}&  0\\
0& 0& 2R_a^{\theta-} & 0
\end{matrix}\right]
\]
as needed
\end{proof}
\section{Spectra the Hermitian Blocks}
As we addressed earlier in remark \ref{rem:BlockStruc}, the quaquaversal operator, $z$, is expressed as
\begin{equation}
1/8\begin{bmatrix}
2R_a^{\pi/2+} + 2R_b^{\pi/3+}+4I& 0& 0& 0 \\
0& -2R_b^{\pi/3+}+2I& 0& 0 \\
2R_b^{\pi/3-}& 0& 2I& 0 \\
0& -2R_b^{\pi/3-}& 2R_a^{\pi/2-} & 4I
\end{bmatrix}
\end{equation}
with respect to 
\[\Pi_{a, b} = \Lambda_{1,1},\Lambda_{1, -1},\Lambda_{-1, 1},\Lambda_{-1, -1}\]
with
\[\Lambda_{\alpha, \beta} = \{v| R_a^\pi v = \alpha v\text{ and } R_b^\pi v = \beta v \}\]

As before, $R_a^{\theta+} = \frac{R_a^\theta+R_a^{-\theta}}{2}$ and $R_a^{\theta-} = \frac{R_a^\theta-R_a^{-\theta}}{2}$. And as before, we will write $S = R_a^{\pi/2}$ and $T = R_b^{\pi/3} $. It is not difficult to the eigenvalues of three of the four blocks along the main diagonal, the $\Lambda_{-1, -1}$ block, the $\Lambda_{1,-1}$ block, and the $\Lambda_{-1, 1}$ block. The last block, $\Lambda_{+1, +1}$ is quite difficult to analyze and contains all the eigenvalues near $1$. We shall now compute the eigenvalues of the three easy blocks.

The easiest blocks are $1/8(2I) = (1/4)I$ on $\Lambda_{-1, 1}$ and $1/8(4I) = (1/2)I$ on $\Lambda_{-1, -1}$. Since these are just multiples of the identity we see that the blocks have eigenvalues $1/4$ with multiplicity equal to the dimension of $\Lambda_{-1, 1}$ and $1/2$ with multiplicity equal to the dimension of $\Lambda_{-1, -1}$, respectively. Since the trace of a projection is equal to the dimension of the image we just need to find projections onto each of these subspaces and compute their traces. Let us verify that 
\[P_{\alpha, \beta} = \frac{1}{4}(I+\alpha R_a^\pi + \beta R_b^\pi + \alpha\beta R_c^\pi)\]
is a projection onto $\Lambda_{\alpha,\beta}$. Here $R_c^\pi =  R_a^\pi R_b^\pi = R_b^\pi R_a^\pi$, is a rotation by $\pi$ about the axis, $c$, perpendicular to $a$ and $b$. 

Let $v\in\Lambda_{\alpha',\beta'}$ with $\alpha', \beta'$ both in $\{\pm 1\}$. Then
\[P_{\alpha, \beta}v = 1/4(1 + \alpha\alpha' + \beta\beta'+ \alpha\beta\alpha'\beta')v\]

If $(\alpha,\beta,\alpha',\beta')\in\{\pm1\}^4$ then $(1 + \alpha\alpha' + \beta\beta'+ \alpha\beta\alpha'\beta')$ vanishes whenever $(\alpha, \beta)\neq (\alpha',\beta')$ and is equal to $4$ otherwise. Therefore, $P_{\alpha, \beta}$ projects onto $\Lambda_{\alpha, \beta}$ as needed.

The trace of $R_b^\pi$ is $(-1)^k$, so we have: 
\[\dim\Lambda_{\alpha, \beta} = \tr P_{\alpha\beta} = 1/4(2k+1+\alpha(-1)^k+\beta(-1)^k+\alpha\beta(-1)^k)\]
and therefore,
\[\dim\Lambda_{1, 1} = \tr P_{11} = 1/4(2k+1+3(-1)^k)\]
and,
\begin{equation}
\dim\Lambda_{-1, 1} = \dim\Lambda_{1, -1} = \dim\Lambda_{-1, -1} = 1/4(2k+1-(-1)^k)
\end{equation}
Since $z$ is $(1/4)I$ on $\Lambda_{-1, 1}$ and $(1/2)I$ on $\Lambda_{-1, -1}$ these blocks have eigenvalues $1/4$ and $1/2$ respectively with multiplicity $1/4(2k+1-(-1)^k)$.

It is a bit more difficult to compute the $\Lambda_{1,-1}$ block. The Quaquaversal operator, $z$, is $1/8(2I - 2R_b^{\pi/3+}) = 1/4(I - R_b^{\pi/3+}) = 1/4(I - T^+)$ on this block.

We will show that the $1/2$-eigenvectors of $T^+$ in the $\Lambda_{1,-1}$ block are exactly those vectors $u=v_\eta +S^2v_\eta$ where $\eta$ is a primitive sixth root of unity and $v\in H_{2k+1}$ satisfies $Tv_\eta = \eta v_\eta$. It will turn out that the only eigenvalues of $T^+$ in the $\Lambda_{1,-1}$ block are $1/2$ and $-1$. Since we already know that the dimension of this block is $1/4(2k+1-(-1)^k)$ computing the multiplicity of the $1/2$ eigenvalue will allow us to determine the multiplicity of $-1$ by subtraction.

Let us now verify that the $v_\eta +S^2v_\eta$ are indeed $1/2$-eigenvectors of $T^+$ in the $\Lambda_{1,-1}$ block.

The $\Lambda_{1, -1}$ block consists of exactly those vectors, $v$, satisfying:
\begin{enumerate}[1.)]
\item $R_b^\pi v = T^3v = -v$
\item $R_a^\pi v= S^2v = v$
\end{enumerate}
Recall that $S^2T = T^{-1}S^2$.
\[T^3(v_\eta + S^2v_\eta) = \eta^3v_\eta + T^3S^2v_\eta= -v_\eta + S^2T^{-3}v_\eta =-v_\eta + S^2\bar{\eta}^3v_\eta = -(v_\eta + S^2v_\eta)\]
\[S^2(v_\eta + S^2v_\eta) = S^2v_\eta + S^4v_\eta= v_\eta + S^2v_\eta\]
We have now verified that $v_\eta +S^2v_\eta$ is in $\Lambda_{1,-1}$. Let us now verify that $T^+(v_\eta +S^2v_\eta) = 1/2(v_\eta +S^2v_\eta)$. Using the fact that $\frac{\eta +\bar{\eta}}{2}=\frac{(1+\sqrt{3}i)/2 + (1-\sqrt{3}i)/2}{2}= 1/2$ we get:
\begin{multline}T^+(v_\eta +S^2v_\eta) =\frac{T+T^{-1}}{2}(v_\eta +S^2v_\eta)=\frac{\eta+\bar{\eta}}{2}v_\eta+
\frac{T+T^{-1}}{2}S^2v_\eta=\\
\frac{1}{2}v_\eta+
S^2\frac{T^{-1}+T}{2}v_\eta =\frac{1}{2}(v_\eta+
S^2v_\eta)
\end{multline}
as needed.

We shall now prove the converse---that the elements of $\Lambda_{1, -1}$ which are 1/2-eigenvectors of $T^+$ are of the form $v_\eta + S^2v_\eta$.

Since $T$ is diagonalizable and $T^3 = -I$ on $\Lambda_{1,-1}$, we know that any $v\in\Lambda_{1,-1}$ can be uniquely written as $v = v_\eta + v_{\bar{\eta}} + v_{-1}$ with $Tv_\gamma = \gamma v_\gamma$. Furthermore, $v_{\bar{\eta}} = S^2v_\eta$. To see this, notice that $S^2$ swaps the $\eta$-eigenspace of $T$ with the $\bar{\eta}$-eigenspace. Indeed, 
\[T(S^2v) = S^2T^{-1}v = S^2\bar{\eta}v = \bar{\eta}(S^2v)\]
So,
\[v_\eta + v_{\bar{\eta}} + v_{-1} = v = S^2v = S^2v_{\bar{\eta}} + S^2v_\eta + S^2v_{-1} \]
Since the representation is unique, we know that $S^2v_\eta = v_{\bar{\eta}}$ as needed.

If $v$ is an eigenvector of $T^+$, then it is either of the form $v_\eta + v_{\bar{\eta}}$ or $v_{-1}$. Indeed, 
\begin{multline}
T^+(v_\eta + v_{\bar{\eta}} - v_{-1}) =  \frac{\eta+\bar{\eta}}{2}(v_\eta + v_{\bar{\eta}}) - v_{-1} =\\
\frac{(1+\sqrt{3}i)/2 + (1-\sqrt{3}i)/2}{2}(v_\eta + v_{\bar{\eta}})-v_{-1} = 1/2(v_\eta + v_{\bar{\eta}}) -v_{-1}
\end{multline}
So to be an eigenvector either $v_\eta$ and $v_{\bar{\eta}}$ are both zero or $v_{-1}$ is zero. In the first case the eigenvalue is $-1$ and in the second case it is $1/2$. 

Combining these two facts we see that the 1/2-eigenvectors of $T^+$ in $\Lambda_{-1, 1}$ are of the form $v_\eta + S^2v_\eta$, as needed.

We have just shown that the space of 1/2-eigenvectors of $T^+$ in $\Lambda_{-1, 1}$ is $I+S^2$ applied to the space of $\eta$-eigenvectors of $T$. As $I+S^2$ is injective, indeed if $v_\eta + S^2v_\eta = v_\eta' + S^2v_\eta'$ then $v_\eta - v_\eta' + S^2v_\eta-S^2v_\eta=0$ implies $v_\eta -v_\eta' =0$ as $v_\eta - v_\eta'$ and $S^2v_\eta-S^2v_\eta$ live in disjoint subspaces. 
On the $\eta$-eigenvectors of $T$, the dimension of this space is the dimension the space of 1/2-eigenvectors of $T^+$ in $\Lambda_{-1, 1}$ is equal to the multilicity of $\eta$-eigenvectors of $T$ in $H_{2k+1}$. $T$ has eigenvalues given by $\{\eta^j\}_{j=-k}^{j=k}$ so the $\eta$ occurs with multiplicity $\lfloor\frac{k+5}{6}\rfloor + \lfloor\frac{k+1}{6}\rfloor$. The remaining eigenvalue of $T^+$ on this block is $-1$, so to compute it's multiplicity we just subtract this from the dimension of the block and get \[1/4\left(2k+1-(-1)^k\right)-\left(\lfloor\frac{k+5}{6}\rfloor + \lfloor\frac{k+1}{6}\rfloor\right)\]
So, since $z=1/4(I-T^+)$ on this block. We have $1/4(1-1/2) = 1/8$ with multiplicity $\lfloor\frac{k+5}{6}\rfloor + \lfloor\frac{k+1}{6}\rfloor$ and $1/4(1-(-1)) = 1/2$ with multiplicity $1/4\left(2k+1-(-1)^k\right)-\left(\lfloor\frac{k+5}{6}\rfloor + \lfloor\frac{k+1}{6}\rfloor\right)$
\bibliographystyle{plain}
\bibliography{references.bib}

\begin{thebibliography}{1}

\bibitem{bourgain2008spectral}
Jean Bourgain and Alex Gamburd.
\newblock On the spectral gap for finitely-generated subgroups of su (2).
\newblock {\em Inventiones mathematicae}, 171(1):83--121, 2008.

\bibitem{bump2004lie}
Daniel Bump.
\newblock {\em Lie groups}, volume~8.
\newblock Springer, 2004.

\bibitem{conway1998quaquaversal}
John~H Conway and Charles Radin.
\newblock Quaquaversal tilings and rotations.
\newblock {\em Inventiones mathematicae}, 132(1):179--188, 1998.

\bibitem{draco2000growth}
Brimstone Draco, Lorenzo Sadun, and Douglas Van~Wieren.
\newblock Growth rates in the quaquaversal tiling.
\newblock {\em Discrete \& Computational Geometry}, 23(3):419--435, 2000.

\bibitem{knapp2016representation}
Anthony~W Knapp.
\newblock {\em Representation theory of semisimple groups}.
\newblock Princeton university press, 2016.

\bibitem{weyl2016classical}
Hermann Weyl.
\newblock {\em The classical groups}.
\newblock Princeton university press, 2016.

\end{thebibliography}

\end{document}